%



\def\ms{\medskip}

\def\phi{\varphi}
\def\cov{{\rm cov}}

\def\su{\subseteq}
\def\a{\alpha}
\def\b{\beta}

\def\d{\delta}
\def\l{\lambda}
\def\k{\kappa}

\def\om{\omega}

\def\lng{\langle}
\def\rng{\rangle}
\def\ov{\overline}
\def\sm{\setminus}
\def\cont{{2^{\aleph_0}}}



\def\cf{{\rm cf\,}}

\def\otp{{\rm otp\,}\,}

\def\id{{\rm id}}

\def\fin{{\rm Fin}}
\def\Fin{\fin}


\def\endproof#1{\hfill {\parfillskip0pt$\smiley_{\hbox{{#1}}}$\par\medbreak}}

\def\imply{\Rightarrow}

\def\proof{\smallbreak\noindent{\sl Proof}: }

\def\Cal#1{{\cal #1}}

\newcount\itemno
\def\itm{\advance\itemno1 \item{(\number\itemno)}}
\def\ritm{\advance\itemno1 \item{)\number\itemno(}}
\def\startitm{\itemno=-1 }
\def\startaitm{\itemno=0}

\def\aitm{\advance\itemno1 
\item{(\letter\itemno)}}

\def\letter#1{\ifcase#1 \or a\or b\or c\or d\or e\or f\or g\or h\or
i\or j\or k\or l\or m\or n\or o\or p\or q\or r\or s\or t\or u\or v\or
w\or x\or y\or z\else\toomanyconditions\fi}

\newcount\secno

\magnification=1200
\baselineskip 18pt
\def\cp{{\rm cp\,}}
\def\Forb{{\rm Forb\,}}
\def\F{{\Cal F}}

\headline{}

\def\rest{|}

\def\wcp{{\rm wcp\, }}

\def\Bbb{\bf}
\def\endproof{\hfill$\triangle$\ms}

\font\bigfont cmbx10 scaled \magstep2
\font\namefont cmbx10 scaled \magstep2
\font\abstfont cmbx10 scaled 800
{{
\obeylines

\everypar={\hskip0cm plus 1 fil}
\bf
\parskip=0.4cm

{\bigfont REPRESENTING EMBEDDABILITY}
{\bigfont AS SET INCLUSION}

\medskip
{\rm June 1995}
\bigskip

\parskip0.1cm

\medskip

{\namefont Menachem Kojman}
Department of Mathematics 
Carnegie-Mellon University
Pittsburgh PA, 15213
{\tt kojman@andrew.cmu.edu}

}
\vfill
\everypar{}
\rm
\leftskip1cm\rightskip1cm

ABSTRACT. {\abstfont A few steps are made towards representation theory
of embeddability among uncountable graphs. A monotone class of graphs
is defined by forbidding countable subgraphs, related to the graph's
end-structure.  Using a combinatorial theorem of Shelah it is proved:

\item{-} The complexity of the class in every regular uncountable
$\l>\aleph_1$ is at least $\l^++\sup\{\mu^{\aleph_0}:\mu^+<\l\}$

\item{-} For all regular uncountable $\l>\aleph_1$ there are $2^\l$
pairwise non embeddable graphs in the class having strong homogeneity
properties.

\item {-} It is characterized when some invariants of a graph $G\in
\Cal G_\l$  have to be inherited by one of fewer than $\l$ subgraphs whose
union covers $G$.

All three results are obtained as corollaries of a {\it representation
theorem} (Theorem 1.10 below), that asserts the existence of a
surjective homomorphism from the relation of embeddability over
isomorphism types of regular cardinality $\l>\aleph_1$ onto set
inclusion over all subsets of reals or cardinality $\l$ or less.
Continuity properties of the homomorphism are used to extend the first
result to all singular cardinals below the first cardinal fixed point
of second order.

The first result shows that, unlike what Shelah showed in the class of
all graphs, the relations of embeddability in this class is
not independent of negations of the GCH. }
\footline{\hfill}
\global\pageno0

\eject
}

\headline{\tensl\hfill M.~ Kojman: embeddability\hfill}

{\bf \S0 Introduction}

The study of embeddability among infinite structures has a long
tradition of invoking combinatorics. One well known example is Laver's
use of Nash-Williams' combinatorial results to show that embeddability
among countable order types is well quasi ordered [L]. In the  study of
embeddability among uncountable structures, the most prominent
combinatorial principle has been the Generalized Continuum Hypothesis (GCH),
which asserts that every infinite set has the least possible number of
subsets.  Hausdorff proved as early as 1914 using the GCH that in every
infinite power there is a universal linear ordering, that is, one in
which every linear ordering is embedded as a subordering.  Jonsson [Jo]
used the GCH to prove that classes of structures satisfying a list of
6  axioms have universal structures in all uncountable powers. See
also [R] for graph theory and [MV] for model theory.


Finer combinatorial principles have come from Jensen's work in
G\"odel's universe of constructible sets [Je]. Thus, for example,
Macintyre [M] uses Jensen's diamond --- a principle stronger than CH
--- to prove that no abelian locally finite group of size $\aleph_1$
is embeddable in all universal locally finite groups of size
$\aleph_1$, and Komjath and Pach [KP1] use the same principles to
prove that there is no universal graph in power $\aleph_1$ among all
graphs omitting $K_{\om,\om_1}$.

A common property of the combinatorial principles mentioned above is that
they are not provable from the usual axioms of Set Theory. 
Easton [E] showed that the GCH can fail for all regular cardinal.
Magidor [Ma] showed the GCH could fail at $\aleph_\om$ with GCH below
it and Foreman and Woodin [FW] showed that the GCH could fail everywhere
(both using large cardinals). In Spite of this, the common impression among
mathematicians working in areas having intimate relations to infinite
cardinals, like infinite graph theory, infinite abelian
groups, and model theory, remained  that the GCH was a useful assumption,
while its negations were not. 

In the context of embeddability this impression was fortified by
Shelah's independence results. Shelah showed that universal structures in
uncountable powers may or may not exist under negations of the GCH.
Thus, while GCH implies the existence of universal graphs in all
infinite $\l$, the assumption $\l<\cont$ for regular uncountable $\l$
does not determine the existence or non existence of a universal graph
in $\l$ (see [S3], [Me] and [K]). 

Shelah's independence results [S1,2,3] created the expectation that
the existence of a universal structure in a class of structures in
uncountable cardinalities would always be independent of negations of
GCH, unless the existence was trivial (because the class of structures
is ``dull''). See, for example, [KS] for results about the class of
$K_{\om_1}$-free graphs that support this expectation.

The understanding of negations of the GCH at singular cardinals has
changed dramatically in the last five years.  The most fascinating
development in this area is Shelah's bound on the exponents of
singular cardinals. Shelah proved the following magnificent theorem,
formulated here, though, in a way Shelah himself resents:

\ms\noindent{\bf 
0.1 Theorem}:
{\sl
 If $2^{\aleph_n}<\aleph_\om$ for all $n$ then
$2^{\aleph_\om}<\aleph_{\om_4}$.
}

\ms
Knowing that by Cohen's results no bound can be put on the exponent of
a regular cardinal, this theorem is exceptionally thrilling. A short
proof of it can be found in [J].

The formulation Shelah prefers is the following:
\ms\noindent{\bf
0.1a Theorem}: $\cf([\aleph_\om]^{\aleph_0},\su)<\aleph_{\om_4}$

\ms
This formulations says that the cofinality of the partial ordering of
set inclusion over countable subsets of $\aleph_\om$ is ALWAYS smaller
than $\aleph_{\om_4}$, no matter how large $\aleph_\om^{\aleph_0}$ may
be. In other words, this theorem exposes a robust {\it structure} of
the partial ordering of set inclusion, which is affected by negations
of the GCH in a limited way only.  The reader will verify that 0.1a
implies 0.1. A proof of this theorem is in Shelah's recent book on
Cardinal arithmetic [S]. In this book Shelah reduces the problem of
computing the exponent of a singular cardinal to an algebra of reduced
products of regular cardinal, and uses a host of new and sophisticated
combinatorics to analyze the structure of such reduced products.

A common property of the combinatorial principles Shelah uses in [S]
and in later works on cardinal arithmetic, is that they {\it are}
proved in ZFC, the usual axiomatic framework of set theory. This is
necessary, since 0.1a (unlike the conclusion of 0.1) is an absolute
theorem, namely assumes nothing about cardinal arithmetic.

In this paper we use some of Shelah's combinatorics to expose robust
connections between the structure of embeddability over a monotone
class of infinite graphs and the relation of set inclusion. This is
done by means of a {\it representation theorem}, that asserts the
existence of a surjective homomorphism from the former relation onto
the latter. One corollary is that the structure of embeddability over
the class we shall study --- which is defined by imposing restrictions
on the the graph's end-structure --- is not independent of negations
of GCH, but also information that is not related to cardinal
arithmetic is obtained.

Shelah's ZFC combinatorics on uncountable cardinals was found useful
in the study of embeddability in several papers. In [KjS1] it was
shown that if $\l>\aleph_1$ is regular and $\l<\cont$ then there is no
universal linear ordering in $\l$. In other words, an appropriate
negation of CH determines negatively the problem of existence of a
universal linear ordering in power $\l$. Similar results were proved
for models of first order theories [KjS2]; infinite abelian groups
[KjS3] and [S4] and metric spaces [S5]. But so far no application was
found for infinite graphs, in spite of the existing rich and active
theory of universal graphs.

The theory of universal graphs, that  began with Rado's construction
[R] of a countable strongly universal graph,  has advanced
considerably since, especially in studying universality over monotone
classes (see [DHV] for motivation for this). A monotone class of
graphs is always of the form $\Forb(\Gamma)$, the class of all graphs
omitting a some class $\Gamma$  of ``forbidden'' configurations as
subgraphs. A good source for the development of this theory is the
survey paper [KP1], in which the authors suggest a generalization of
universality, which they name ``complexity'': the complexity of a
class of graphs is the least number of members in the class needed to
embed as induced subgraphs all members in the class. The complexity is
1 exactly when a strongly universal graphs exists in the class.

\ms
The paper is organized as follows. In Section 1 a class of graphs is
specified by forbidding countable configurations related to the
graph's end-structure, and it is noted that by a generalization of a
theorem by Diestel, Halin and Vogler the complexity of the resulting
class $\Cal G$ at power $\l$ is at least $\l^+$. A surjective
homomorphism is now constructed from the relation of (weak)
embeddability over $\Cal G_\l$ for regular $\l>\aleph_2$ onto the
relation of set-inclusion over all subsets of reals of cardinality
$\le
\l$. Combining both results, $\max\{\l^+,\cont\}$ is set as a lower
bound for the complexity of $\Cal G_\l$ for regular $\l>\aleph_1$.

In Section 2 a certain continuity property of the homomorphism from
Section 1 is proved, and is used to extend the lower bound from
Section 1 to all singular cardinals below the first fixed point of
second order. In this Section the representation Theorem is stated in
its full generality, generalizing Theorem 1.8 to higher
cardinals.

In Section 3 it is proved that in every regular $\l>\aleph_1$ there
are $2^\l$ pairwise non mutually embeddable elements in $\Cal G_\l$,
each of which being ``small'' in the sense that it is mapped by the
homomorphism to a finite set.  For the case $\l$ inaccessible, this
result makes use of a very recent result by Gitik and Shelah about
non-saturation of the non-stationary ideal on $\l$. No cardinal
arithmetic assumptions are made in this Section and in Section 4.

In Section 4 a decomposition theorem is proved for a proper subclass
of $\Cal G_\l$, $\l>\aleph_1$ regular, which is also defined by
forbidding countable configurations. The Theorem gives a necessary and
sufficient condition to when the invariant of a graph $G$ in the class
is inherited by at least one subgraphs from a collection of $<\l$
subgraphs whose union covers $G$.

\bigskip
\noindent
{\bf NOTATION} 
A graph $G$ is a pair $\lng V,E\rng$ where $V$ is the set of vertices
and $E\su [V]^2$ is the set of edges. By $G[v]$ we denote the {\it
neighbourhood} of $v\in V$ in $G$, namely $\{u\in V: \{v,u\}\in E\}$.
A graph $G$ is {\it bipartite} if there is a partition $G=G_1\cup G_2$
of $G$ to two (non-empty) disjoint independent vertex sets, each of which is
called a {\it side}.

An ordinal is a set which is well ordered by $\in$. A {\it cardinal}
is an initial ordinal number. The cardinality $|A|$ of a set $A$ is
the unique cardinal equinumerous with $A$. The {\it cofinality} $\cf \l$ of a
cardinal $\l$ is the least cardinal $\k$ such that $\l$ can be
represented as a union of $\k$ sets, each of cardinality less than
$\l$. A cardinal $\l$ is singular if $\cf \l<\l$ and is {\it regular}
if $\cf \l=\l$. If $\k,\k'$ are cardinals we denote by $K_{\k}$
the complete graphs on $\k$ vertices and by $K_{\k,k'}$ the complete
bipartite graphs with $\k$ vertices in one  side and $\k'$ in the other.

If $G_1$ is isomorphic to a subgraph of $G_2$ we write $G_1\le_w G_2$
and we write $G_1\le G_2$ if $G_1$ is isomorphic to an {\it induced}
subgraph of $G_2$. We also say the $G_1$ is {\it embeddable}
(embeddable as an induced subgraph) if $G_1\le_wG_2$ ($G_1\le G_2$).

Classes of graphs will be denoted by $\Cal G$ and $\Gamma$ and are
always assumed to be closed under isomorphism.  A class of graphs
$\Cal G$ is {\it monotone} if $G_1\le_w G_2\in \Cal G \imply G_1\in
\Cal G$. If $\Gamma$ is a set of graphs then $\Forb(\Gamma)$ is the
class of all graphs without a subgraph in $\Gamma$.  Let $\cal G_\l$
be the set of all isomorphism types of $\Cal G$ whose cardinality is
$\l$.  The relations $\le$ and $\le_w$ are reflexive and transitive,
and therefore $\lng \Cal G_\l,\le\rng$ and $\lng \Cal G_\l,\le_w\rng$
are quasi-ordered sets for all classes $\Cal G$ and cardinals $\l$.

Let $\cp\Cal G_\l$, the {\it complexity} of $\Cal G_\l$, be the least
cardinality of a subset $D\su \Cal G_\l$ with the property that for
every $G\in \Cal G_\l$ there exists $G'\in D$ such that $G\le G'$; the
{\it weak complexity} is defined by replacing $\le$ by $\le_w$ (see
[KP]). Clearly, $\wcp \Cal G\le \cp \Cal G$ for any class $\Cal
G$.

The complexity $\cp \Cal G_\l$ is 1 iff there is a graph $G^*\in \Cal
G_\l$ with the property that every member of $\Cal G_\l$ is isomorphic
to an induced subgraph of $G^*$. Such a graph $G^*$ is called {\it
universal in $\l$} (or, sometimes, ``strongly universal'') for the
class $\Cal G$. $\wcp \Cal G_\l=1$ is equivalent to the existence of a
{\it weakly universal} element in $\Cal G_\l$.

Suppose that $A$ is a given infinite set. By
$[A]^\l$ we denote the collection $\{B:B\su A\;\&\;|B|=\l\}$ of all
subsets of $A$ whose cardinality is $\l$. If $B_1\in [A]^\l$ is
contained as a subset in $B_2\in [A]^\l$ we write $B_1\su B_2$. Since
$\su$ is reflexive, transitive and antisymmetric, $\lng [A],\su\rng$
is a partially ordered set. Let $\cov  A_\l$, the covering number of
$[A]^\l$, be the least cardinality of a subset $D\su [A]^\l$
with the property that for every $B\in [A]^\l$ there exists $B'\in D$
such that $B\su B'$.

We remark that the least cardinality of a dominating subset is defined
for every quasi-ordered set, and bears the name ``cofinality''; but we
stick here to the customary graph-theoretic and set theoretic existing
terminologies and refer to the former as ``complexity'' and to the
latter as ``covering number''.

\ms\noindent{\bf 
0.2 Definition}: Let $\l$ be an uncountable regular
cardinal. A {\it club} of $\l$ is a closed (in the order topology)
and unbounded subset of $\l$. Club sets generate a {\it filter} over
$\l$, indeed a {\it $\l$-complete} filter: the intersection of fewer
than $\l$ subsets of $\l$, each of which contains a club, contains a
club. A subset of $\l$ is called {\it stationary} if its intersection
with every club of $\l$ is non empty. The ideal of all subsets of $\l$
which are disjoint to some club of $\l$ is the {\it non-stationary}
ideal. Thus club sets are analogous to measure 1 sets, non-stationary
sets are measure zero and stationary sets are positive measure
(meet every measure 1 set). Let $S^\l_\k$ be $\{\a:\a<\l\wedge
\cf\a=\k\}$ and $S^\l_0=\{\a:\a<\l\wedge\cf\a=\om\}$. 

We shall  need the following combinatorial tool:

\ms\noindent{\bf 
0.3 Theorem}:(Shelah)
{\sl If $\l>$ is  regular, $\mu$ a
cardinal and $\mu^+<\l$  then there is a stationary set $S\su
\l$ and a sequence $\bar C=\lng c_\d:\d\in S\rng$ with $\otp c_\d=\mu$
and $\sup c_\d=\d$ such that for every closed unbounded $E\su
\l$ the set $N(E):=\{\d\in S: c_\d\su E\}$ is stationary. 
}

For a proof see  I[Sh-e, new VI\S2]= [Sh-e, old III\S7].

A sequence $\bar C$ as in the theorem is called a ``club guessing
sequence''. If the $c_\d$ are thought of as ``guesses'', then the
theorem says that for every club (measure 1) set stationarily many
(positive measure) of the guesses are successful.

Suppose $\bar C$ is a club guessing sequence as above.  We define two
{\it guessing ideals} over $\l$, $\id^a(\bar C)$ and $\id^b(\bar C)$,
as follows: 

\ms\noindent{\bf 
0.4  Definition}: 
\startitm
\itm $X\in \id^a (\ov C)$ iff for some club $E\su
\l$ it holds that $c_\d\not\su E$ for all $\d\in S\cap E$. 

\itm $X\in \id^b (\ov C)$ iff for some club $E\su
\l$ it holds that $c_\d\not\su^* E$ for all $\d\in S\cap E$, where
$c_\d\su^*$ means that an end segment of $c_\d$ is contained in $E$.

Thus a set $X\su \l$ is in $\id^a(\bar C)$ iff there are no
stationarily many $\d\in X$ such that $c_\d$ is contained in $E$ for
some club $E$, and $X\su \l$ is in $\id^b(\bar C)$ iff there are no
stationarily many $\d\in X$ such that $c_\d$ is almost (=except for a
proper initial segment) contained in $E$ for some club $E$.

 The ideal $\id^a(\bar C)$ is a $\l$-complete ideal over $\l$ and
$\id^b(\bar C)$ is normal. Also, $\id^b(\bar C)\su\id^b(\bar C)$.

\ms\noindent{\bf 
0.5  Definition}: Let $\om$ be the set of natural numbers. Let
$\Fin$ be the set of all finite subsets of $\om$. Two subsets $X,Y\su
\om$ are {\it equivalent} mod $\Fin$ iff the symmetric difference $X\sm
Y\cup Y\sm X\in \Fin$. By $\Cal P(\om)$ we denote the {\it power set}
of $\omega$ and by $\bar {\Cal P}(\om)$ we denote ${\Cal P}(\om)/\Fin$ the
set of all equivalence classes of subsets of $\om$ modulo $\Fin$. 

Finally, we need a few definitions about reduced powers. A reduced
power is a generalization of {\it ultra-power}. 

\ms\noindent{\bf 
0.6 Definition}: Suppose that $A$ is a structure, $\l$ a cardinal and
$F$ a filter over $\l$. Let $A^\l$ be the set of all functions from
$\l$ to the structure $A$ and let $A^\l/F$ be the reduced power of $A$
modulo $F$.

\bigbreak
\bigbreak
\noindent
{\bf \S1 Representing embeddability as set inclusion}

In [K] it is proved:

\ms\noindent{\bf 
1.1 Theorem}: 
{\sl
If $\Cal G$ is a class of graphs that
contains all $K_{\om,\om}$-free incidence graphs of $A\su \Cal P(\om)$
and the cofinality of the continuum is $\aleph_1$, then $\cp \Cal
G_\l>\cont$ for all uncountable $\l<\cont$.
}

In particular, if $\cf\cont=\aleph_1$ there is no universal graph (in
the class of all graphs) in any uncountable $\l<\cont$.

On the other hand, Shelah proved in [S3]:

\ms\noindent{\bf 
1.2 Theorem}: 
{\sl 
If $\l$ is regular uncountable, it is
consistent that $\l<\cont$ and that a universal graph in power $\l$
exists. 
}
\ms
Mekler [Me] generalized Shelah's result to more general classes of
structures.

Both result together can be understood as follows: a singular $\cont$ affects
the structure of embeddability in a broad spectrum of classes of
infinite graphs below the continuum, but a large regular
$\cont$ may have no effect  on the class of all graphs and
the classes handled by Mekler and Shelah.

It is reasonable to ask if for {\it some} ``reasonably defined'' class
of graphs for which the structure of embeddabilty below $\cont$ is
influenced by the size of $\cont$.  In this section we show that
forbidding certain countable configurations gives rise to a class with
such a desired connection. The configurations we forbid are related to
the {\it end structure} of graphs.

\ms\noindent{\bf 
1.3 Definition}: A {\it ray} in a graph $G$ is a 1-way
infinite path. A {\it tail} of a ray $R\su G$ is an infinite connected
subgraph of $R$. Two rays in $G$ are {\it tail-equivalent} iff they
share a common tail. Tail-equivalence is an equivalence relation on
rays.

We mention in passing that tail-equivalence is a refinement of
{\it end-equivalence}. For more on both relations see [D]. 

\ms\noindent{\bf 
1.4 Definition}: Let $\Cal G$ be the class of all
graphs $G$ satisfying that for every $v\in G$ the induced subgraph of
$G$ spanned by $G[v]$ has at most one ray up to tail-equivalence.

\ms\noindent{\bf 
1.5  Claim}: 
{\sl
There is a non-empty set $\Gamma$ of countable graphs,
each containing an infinite path, such that $\Cal G=\Forb(\Gamma)$. 
}

\proof: 
Let $\Gamma'$ be the set of all countable graphs that contain (at
least) two  tail-inequivalent rays and let $\Gamma$ be all graphs
obtained by choosing an element from $\Gamma'$ and joining a new
vertex to all its vertices. If a graph $G$ contains a subgraph in $\Gamma$
then $G\notin \Cal G$. Conversely, Suppose that $G\notin
\Cal G$. Let $v\in G$ be a vertex such that there are two rays
$R_1,R_2\su G[v]$ which are not tail-equivalent. Let $G'\su G$ be the
induced subgraph spanned by $\{v\}\cup R_1\cup R_2$. Now $G'\in
\Gamma$ and so $G\notin\Forb (\Gamma)$.

Graphs with forbidden countable configuration that contain an infinite
path were considered by Diestel, Halin and Vogler in [DHV] for
$\l=\aleph_0$. They prove (Theorem 4.1):

\ms\noindent{\bf 
1.6 Theorem}: (Diestel-Halin-Vogler)
{\sl
  Let $\Gamma$ be a non-empty set of countable
graphs each containing an infinite path. Then $\Cal G_{\aleph_0}=\Cal
G_{\aleph_0}(\Gamma)$ has no universal element.
}

This Theorem applies to our class $\Cal G$ by 1.5 above, and because every
forbidden configuration in $\Gamma$ contains an infinite path.  The
following is a straightforward generalization of  Theorem 1.6, and
is included only for completeness of presentation's sake:

\ms\noindent{\bf 
1.7  Theorem}:
{\sl
 Let $\Gamma$ be a non-empty set of
countable graphs each containing an infinite path. Then
$\wcp\Forb_\l(\Gamma)\ge\l^+$ for all infinite cardinals $\l$.
}

\proof By induction on $\a<\l^+$ define graphs $G_\a$ as follows:
$G_\a$ is obtained by joining a vertex $w_\a$ to a disjoint union of
 $G_\b$ for all $\b<\a$ (such that $w_\a$ is not in this
union). For all $\l\le\a<\l^+$ the graph $G_\a$ contains no infinite
path and therefore belongs to $\Forb_\l(\Gamma)$. Suppose that $\Cal
F$ is a collection of $\l$ graphs and $f_\a:G_\a\to G(\a)$ is an
embedding of $G_\a$ into some graph $G(\a)\in \Cal F$. By the pigeon
hole principle there is a fixed graph $G\in \Cal F$ such that
$G=G(\a)$ for $\l^+$ many $\a<\l^+$. A second use of the pigeon hole
principle gives  a vertex $w(0)\in G$ and an unbounded set $X\su
\l^+$ such that $w(0)=f_\a(w_\a)$ for all $\a\in X$. This implies, by
the construction of the $G_\a$'s,  that
$G[w(0)]$ contains as subgraphs copies of $G_\a$ for unboundedly many
$\a<\l^+$ and therefore of all  $\a<\l^+$. Repeating this
argument a set $\{w(n):n<\om\}\su G$ is found that spans in $G$ a copy
of $K_{\om}$. Therefore $G$ contains all  countable configurations
and therefore does not belong to $\Forb_\l(\Gamma)$.\endproof

\ms
Thus for every infinite $\l$ we have $\wcp \Cal G_\l\ge \l^+$. Also
Theorem 1.1 from the previous section applies to $\Cal G$,
because $\Cal G$ contains all bipartite graphs; thus (setting
$\theta=\aleph_0$), $\cp \Cal G_\l\ge \cont^+$ if $\cf\cont\le \l$ (we
use here $\wcp\le \cp$).

The virtue of $\Cal G$ is, nevertheless, that $\wcp \Cal G_\l\ge
\max \{\l^+,\cont\}$ regardless to the cardinality of the continuum
for all regular $\l>\aleph_1$ (and many singular $\l$, as seen in the
next section). This is a
corollary of the following:

\ms\noindent{\bf 
1.8 Theorem}: 
{\sl
If $\l>\aleph_1$ is regular then there
is a surjective homomorphism $\Phi:\lng \Cal G_\l,\le_w\rng \to ([\Bbb
R]^{\le \l},\su)$ from the relation of embeddability over
$\Cal G_\l$ onto subsets of reals or cardinality at most $\l$
partially ordered by inclusion.
}

Thus the relation of embeddability among members of $\Cal G_\l$ is at
least as complicated as inclusion among subsets of reals of
cardinality at most $\l$. 

\ms\noindent{\bf 
1.9  Corollary}: 
{\sl
If $\l>\aleph_1$ is regular then $\wcp
\Cal G_\l\ge \max \{\l^+,\cont\}$.
}
\ms
This theorem will be extended to singular values of $\l$ in the next
section.

We turn now to the proof of the theorem. The homomorphism $\Phi$ will
be factored through a reduced product of the inclusion relation over
subsets of reals. We will prove the following  stronger formulation:

\ms\noindent{\bf 
1.10 Theorem}: 
{\sl
Suppose that $\l>\aleph_1$ is
regular. Then  
\startitm
\itm there is a surjective homomorphism $\Phi:\lng \Cal
G_\l,\le_w\rng\to ([\Bbb R]^{\le \l},\su)$

\itm $\Phi$ from $(1)$ can be chosen to be a composition $\psi\phi$
where $\phi$ is a surjective homomorphism to a reduced power
$\left([\Bbb R]^{\le \l},\su\right)^\l/I$ for some normal ideal $I$
over $\l$.
}

\proof: First let us notice that $([\Bbb R]^{\le \l},\su)$ is a
homomorphic image of $\left([\Bbb R]^{\le
\l},\su\right)^\l/I$ for every ideal $I$: Suppose that $\bar A$ is a
representative of an equivalence class of $([\Bbb R]^{\le
\l})^\l/I$. Define $\psi([\bar A]):=\{x\in
\Bbb R: \{\d<\l:x\in A(\d)\}\notin I\}$. In words,
$\psi([\bar A])$ is the set of all reals that appear in a positive set
of coordinates. It is routine to check that the definition of $\psi$
does not depend on the choice of a representative and that $\psi$ is a
homomorphism. 

Thus it suffices to prove that there is a surjective homomorphism
$\phi:\lng \Cal G_\l,\le_w\rng\to \left([\Bbb R]^{\le
\l},\su\right)^\l/I$ for some normal ideal $I$ over $\l$. This
is in fact more than needed for (0). The set $\Bbb R$ can be replaced
here by any set of equal cardinality. It is convenient for us to work
with $\bar {\Cal P}(\om)=\Cal P(\om)/\Fin$.

We shall define a mapping $\phi:\left( \Cal G_\l,\le\right)\to
\left([\bar {\Cal P}(\om)]^{\le \l},\su\right)^\l/I$ after specifying
$I$.  We shall show that $\phi$ is well defined, is a homomorphism and
is surjective. For the definition of the mapping we fix a club
guessing sequence $\ov C=\lng c_\d:\d\in S\rng$, $S\su\l$ stationary
and $\otp c_\d=\om$ for $c_\d$ in $\bar C$. Now let $I=\id^b(\ov C)$.
For each $\d\in S=S^\l_0$ let $\lng \a_n^\d:n<\om\rng$ be the
increasing enumeration of $c_\d$.

Given $G\in \Cal G_\l$ we define $\phi(G)$ after choosing two
auxiliary parameters on $G$. First, we pick a well ordering $<$ of $G$
of order type $\l$, namely a bijection $h$ between the vertices of $G$ and
the ordinals below $\l$, and second, we fix a mapping $r$ so that
$r(v)=\emptyset$ for $v\in G$ in case there are no rays in $G[v]$ and
$r(v)$ is a ray in $G[v]$ otherwise.

Let $G^{\a,<}=\{v\in G: h(v)<\a\}$ for $\a<\l$. When $<$ is unambiguous
we write $G^\a$ for $G^{\a,<}$. For a vertex $v\in G$ and
an ordinal $\d\in S$ we define:

$$\phi_{<,r}(v,\d)=[\{n<\om: r(v)\cap G^{\a_n}\su r(v)\cap
G^{\a_{n+1}}\}]_{\Fin}\leqno{(1)}$$

Thus $\phi_{<,r}(v,\d)$ belongs to $\bar {\Cal P}(\om)$. The
definition in (1) depends strongly on the choice of the well ordering
$<$. In fact, for a vertex $v$ whose neighbourhood $G[v]$ does contain
a ray, $\phi_{<,r}(v,\d)$ can be made any prescribed element of
$\bar{\Cal P}(\om)$ by a suitable choice of $<$. But replacing $r(v)$
by a tail-equivalent $r'(v)$ produces at most a finite change in
$\{n<\om: r(v)\cap G^{\a_n}\su r(v)\cap G^{\a_{n+1}}\}$ and therefore
does not change the definition (1).

Now let 

$$\phi_{<,r}(G,\d)=\{ \phi_{<,r}(v,\d): v\in G\}\leqno{(2)}$$

 Since
$|G|=\l$, we conclude that   

$$\phi_{<,r}(G,\d)\in [\bar {\Cal P}(\om)]^{\le\l}\leqno{(3)}$$

Finally, let 

$$\phi(G)=\phi_{<,r}(G)=[\lng \phi_{<,r}(G,\d):\d\in S\rng]_I\leqno{(4)}$$

The sequence $\lng \phi_{<,r}(G,\d):\d\in S\rng$ belongs to $\left(
[\bar {\Cal P}(\om)]^{\le\l}\right)^\l$, and we let $\phi_{<,r}(G)$ be
the equivalence class of this sequence in the reduced power mod $I$.

\ms\noindent{\bf 
1.11 Proposition}:
\startaitm
\aitm The definition of $\phi(G)$ does not depend on the choice of $<$
and $r$.

\aitm $\phi$ is a homomorphism: if $G_1\le G_2$ then
$\phi(G_1)\su_I\phi(G_2)$.

\aitm $\phi$ is surjective.

\proof: We shall prove (a) and (b) simultaneously by proving 
\aitm If
$G_1\le G_2$ are in $\Cal G_\l$ and $<_1,<_2,r_1,r_2$ are any
parameters as in (1) above for $G_1,G_2$ respectively,  then
$\phi_{<_1,r_1}(G_1)\su_I\phi_{<_2,r_2}(G_2)$

Then (a) follows from (d) 
by putting $G_1=G_2=G$ and (b) follows from (a) and (d).

Let $<_1,<_2,r_1,r_2$ be given. We need to show that
$\phi_{<_1,r_1}(G_1)\su_I\phi_{<_2,r_2}(G_2)$, namely, that:

$$S\sm \{ \d\in S: \phi_{<_1,r_1}(G_1,\d)\su
\phi_{<_2,r_2}(G_2,\d)\}\in I\leqno{(5)}$$

Fixing an embedding from $G_1$ to $G_2$ we assume, without loss of
generality, that $G_1$ is a subgraph of $G_2$. The following set is
closed unbounded in $\l$ by a standard back and forth argument:

$$E=\{\a<\l: G_2^\a\cap G_1=G_1^\a\}\leqno{(6)}$$

Therefore, by the definition of $I$, the set $N(E)=\{\d\in S:c_\d\su
E\}$ satisfies

$$S\sm N(E)\in I\leqno{(7)}$$

We shall show that for every $\d\in N(E)$ we have for all $v\in G_1$
such that $r_1(v)\not=\emptyset$:

$$\phi_{<_1,r_1}(v,\d)=\phi_{<_2,r_2}(v,\d)\leqno{(8)}$$

Let $v\in G_1$ be given. Clearly, $G_1[v]\su G_2[v]$. Thus $r_1(x)$
and $r_2(x)$ are both rays in $G_2(x)$. Since $G_2\in \Cal G$, the
rays $r_1(v),r_2(v)$ are tail-equivalent. Fix a common tail $r(v)$ of
$r_1(v),r_2(v)$. Using $r(v)$ instead of either $r_1(v)$ or $r_2(v)$
in the definition (1) for $\phi_{<_i,r_i}(v,\d)$ ($i=1,2$) makes no
difference. So we may assume without loss of generality that
$r_1(v)=r_2(v)=r(v)$ for all $v\in G_1$.

Suppose that $v\in G_1$ and that $c_\d\su^* E$. Then the sets
$\{n<\om: r(v)\cap G^{\a_n}_1\su r(v)\cap G^{\a_{n+1}}\} $ and
$\{n<\om: r(v)\cap G^{\a_n}_2\su r(v)\cap G^{\a_{n+1}}\}$ are
equivalent modulo finite, because $G_1^{\a_n}=G_2^{\a_n}$ for an end
segment of $\om$. Thus, the definition (2) above gives for every
$\d\in N(E)$:

$$\phi_{<_1,r_1}(G_1,\d)=\{\phi_{<_1,r_1}(v,\d):v\in
G_1\}=\{\phi_{<_2,r_2}(v,\d):v\in G_1\}\su\leqno{(9)}$$

\vskip-30pt
$$\su \{\phi_{<_2,r_2}(v,\d):v\in
G_2\}=\phi_{<_2,r_2}(G_2,\d)$$

And now (5) follows by (7). 

This proves (d) and hence (a) and (b). 

To prove (c) fix a member $\bar A=\lng A_\a:\a<\l\rng \in
\left([\bar{\Cal P}(\om)]^{\le \l}\right)^\l$. We shall construct a
graph $G\in \Cal G_\l$ so that $\phi(G)=[\bar A]_I$.

Enumerate $A_\d=\lng X_{\d,\a}:\a<\l\rng$ for $\d\in S$. By induction
on $\a<\l$ we define an increasing and continuous union of graphs
$G_\a$ such that:
\startaitm
\aitm $|G_{\a}|<\l$
\aitm $\b<\a<\l$ implies that $G_\b$ is an induced subgraph of $G_\a$
and if $\a$ is a limit ordinal then $G_\a=\bigcup_{\b<\a}G_\b$.
\aitm For every $\gamma<\b+1\le \a$ and a vertex $v\in G_\gamma$ there is
a vertex $u\in G_{\b+1}\sm G_\b$ such that
$G_{\b+1}[u]=\{v\}$ ($=G_\a[v]\cap G_{\b+1}$ because $G_{\b+1}$ is an
induced subgraphs of $G_\a$ by (b)).
\aitm if $\a=\d+1$ and $\d\in S$ then for all $\d'\in S\cap (\d+1)$ and
$\gamma\le
\d$ there is a vertex $y(\d',\gamma)\in G_\a\sm G_\d$ and a path
 $\lng z_n:n<\om\rng$ in $G_\d$ such that:
\item\item{(i)} $G_\a[y(\d',\gamma)]=\{z_n:n<\om\}$ and $G_\a[z_n]$
contains no rays.
\item\item{(ii)} If the increasing enumeration of $c_{\d'}\sm
X_{\d',\gamma}$ is $\lng\a_{m(n)}:m<\om\rng$ then $z_n\in
G_{\a_{m(n)+1}}\sm G_{\a_{m(n)}}$. 
\aitm $G_\a\in \Cal G$ and for every $\b<\a$ and $v\in G_{\b+1}\sm
G_\b$ the set $G_\a[v]\sm G_\b$ is independent.

\ms

Suppose first that this construction is carried out, and let
$G=\bigcup_{\a<\l}G_\a$. This is a graph of cardinality $\l$.  By (e)
it follows that $G$ belongs to $\Cal G_\l$: if $v\in G\sm G_0$ let
$\b+1$ be the minimal so that $v\in G_{\b+1}$; $G[v]=G_\b[v]\cup
G[v]\sm G_\b$.  $G_\b[v]$ contains at most one ray up to
tail-equivalence because $G_{\b+1}\in \Cal G$ (condition (e)) and
$G[v]\sm G_\b=\bigcup_{\b<\a<\l}G_\a[v]\sm G_\b$ is an independent set
by (e) and therefore contains no rays at all. If $u\in G_0$ then
$G[u]$ is independent.

 Fix any well ordering $<$ of $G$ of order type $\l$. It is standard
to check that for a closed unbounded set $E\su\l$ it holds that
$G_\a=G^{\a,<}$ for all $\a\in E$. We may restrict attention to this
set of indices alone. Suppose now that $\d'\in S\cap E$ and let
$\gamma<\l$ be given. Find $\d\in S$ so that $\d>\max\{\d',\gamma\}$.
At stage $\d+1$ the vertex $y(\d,\gamma)$ mentioned in (d) satisfies
that $\phi_{<,r}(y(\d,\gamma))=X_{\d,\gamma}$ for all $\gamma<\d$,
where $r$ is any function as in (1) above. This shows that
$\phi_{<,r}(G)=[\bar A]_I$.

Let us see that the induction can be carried out. For $\a=0$ let
$G_\a$ be a single vertex and at limits take unions. Let us check that
conditions (a)--(e) hold for $\a=0$ and hold for a limit $\a$ if they
hold for all $\b<\a$. At successor stages $\a+1$ we distinguish two
cases:

\noindent Case 1: $\a\notin S$. In this case
$G_{\a+1}$ is obtained from $G_\a$ by adding, for each $v\in G_a$, new
vertices $\{x_{v,\b}:\b<\a\}$ and adjoining each of them to $v$; thus
$G_{\a+1}[x_{v,\b}]=\{v\}$. Since we added $|\a|<\l$ new vertices (a)
holds. (b) holds trivially. $G_{\a+1}$ was defined so that (c) holds
and (d) holds vacuously. (e) holds as no two new vertices are joined
by an edge.

\noindent Case 2: $\a=\d\in S$.  For every $v\in G_\a$ add new
$\{x_{v,\b}:\b<\a\}$ exactly as in case 1. In addition add vertices
$y(\d',\b)$ for $\d'\le \d$ and $\b\le \a$.  We specify the neighbours
of $y(\d',\b)$: Let $\lng\gamma_n:n<\om\rng\su C_{\d'}$ be the
increasing enumeration of $X_{\d,\b}$.  By induction on $n<\om$ choose
vertices $z_n\in G_{\gamma_n}$ so that $z_{n+1}$ is connected to $z_n$
and $G_\d[z_n]$ contains no rays. This is possible by condition (b)
and the induction hypothesis, that implies that only the vertices
$y(\d',\b)$ contain rays in their neighbourhoods.  Then connect $y(\d',\b)$ to all $z_n$. The requirement that
$z_n$ has no rays in its neighbourhood is not needed before Section 4,
where it is needed to show that the graphs constructed here lies in a
proper subclass. For the purpose of this proof it can be
ignored.  Conditions
(a)--(c),(e) hold as in the previous case. Condition (d) was just
handled. \endproof

\bigbreak
\bigbreak
\noindent{\bf \S2 Continuity, singulars and $r$-subgraphs}

In this section we study the homomorphism $\Phi$ from Theorem
1.8 and show it has a certain continuity property. 

As a result we will be able to prove that Theorem 1.9 holds
also for many singular cardinals $\mu$. 
\ms\noindent{\bf 
2.1 Definition}: For $G_1, G_2\in\Cal G_\l$ say that $G_1$ is an
$r$-subgraph of $G_2$ iff $G_1$ is a subgraph of $G_2$ and for all
$v\in G_1$ if $G_2[v]$ contains a ray then $G_1[v]$ contains a ray.
Equivalently, if $r$ is a function on $G_2$ such that $r(x)$ is a ray
in $G_2[x]$ if such a ray exists and $\emptyset$ otherwise, then
$\forall x\in G_1 \left (r(x)\su^* G_1\right )$.

\ms\noindent{\bf 
2.2 Claim}: (Continuity)  
{\sl
If $G\in \Cal G_\l$ and
$G=\bigcup_{\a<\b}G_\a$ for some $\b<\l$ so that  $G_\a$ is an
$r$-subgraph of $G$
for every $\a<\b$ then

$$\Phi(G)=\bigcup_{\a<\b} \Phi(G_\a)$$
}

\proof For every $\a<\b$ the relation $\Phi(G_\a)\su_I\Phi(G)$ follows
 because $\Phi$ is a homomorphism and  $G_\a\le G$. 

On the other hand, suppose that $A\in \Phi(G)$, and we will show that
for some $\a<\b$ we have $A\in \Phi(G_\a)$. 

Let $<,r$ and $<_\a,r_\a$
be chosen parameters for $G$ and for  every $G_\a$ for $\a<\b$
respectively, as in the definition of $\phi$. For every $\a<\b$ there
is a set $X_\a\in I$ such that for every $\d\in \l\sm X_\a$ and every
$v\in G_\a$ for which $r(v)\not=\emptyset$ condition (8) in the proof
of Proposition 1.11 holds:

$$\phi_{<_\a,r_\a}(v,\d)=\phi_{<,r}(v,\d)\leqno{(8)}$$

Let $X=\bigcup_{\a<\b} X_\a$. By $\l$-completeness of $I$ we know that
$X\in I$. Suppose $\d\in \l\sm X$ let $v\in G$.

Since $A\in \Phi(G)$, the set $Y:=\{\d<\l: A\in \phi_{<,r}(G,\d)\}$ is
positive. For every $\d\in Y$ pick $v_\d\in G$ so that
$\phi_{<,r}(v_\d,\d)=A$. 

For every $\d\in Y$ there is some $\a(\d)<\b$ for which $v_\d\in
G_{\a(\d)}$. Since $I$ is $\l$-complete, there is a fixed $\a<\b$ and
a positive set $Y'\su Y$ such that $\a(\d)=\a$ for all $\d\in Y'$. 

Since $X\in I$, we may, without loss of
generality, subtract this union from $Y'$, and $Y'$ would still be
positive. But now it follows from (8) that for every $\d\in Y'$ it holds that

$$\phi_{<_\a,r_\a}(v_\d,\d)=\phi_{<,r}(v_\d,r)=A$$

Since $Y'$ is positive this shows that $A\in \Phi(G_\a)$ and thus
completes the proof.

\ms\noindent{\bf 
2.3 Corollary}: 
{\sl
If  $\mu$ is singular and $\mu$ 
smaller than the first cardinal fixed point of second order, then $\cp
\Cal G_\mu\ge\max\{\mu^+,\cont\}$.
}

\proof We address only the term $\cont$, the other following from
1.7. Assume, then, that $\mu$ is singular,  smaller than the
first fixed point of second order, and that $\mu<\cont$. To prove that
$\wcp \Cal G_\mu\ge\cont$ suppose $\k<\cont$ and that for every
$\a<\k$ we are given $G_\a\in \Cal G_\mu$, and we
will exhibit a graph $G\in \Cal G_{\l}$ for some regular $\l<\mu$
which is not embedded in any $G_\a$ for $\a<\k$.

By the assumption on $\mu$, and since $|G_\a|=\mu$ for $\a<\k$, for
every $\a<\cont$ there exists a regular $\l<\mu$ and a family $\Cal
F_\a\su [G_\a]^\l$ with $|\Cal F_\a|=\mu$ such that for all $X\in
[G_\a]^\l$ there is $F\in [\Cal F_\a]^{<\l}$ so that $X\su
\bigcup F$ (see [KjS1],4.5). So while $\F_\a$ itself may not be
dominating in $\lng [G_\a]^\l,\su\rng$, the set of all unions of $<\l$
members of $G_\a$ is dominating (but has cardinality larger than
$\mu$ in this case). We may assume, by increasing each $A\in \Cal
F_\a$ to a larger set of the same cardinality, that each $A\in \Cal
F_\a$ spans an $r$-subgraph of $G_\a$, and abusing notation we shall
not distinguish between $A$ and the subgraph it spans. 

The cardinality of the following set is $\k\times\mu$:

$$D:=\{\Phi(A):A\in\bigcup_{\a<\cont}\F_\a\}$$

As $\k\times\mu<\cont$ and $|\Phi(A)|=\l$ for $\Phi(A)\in D$, we have
$|\bigcup D|<\cont$. Thus we can find a set $B\in [R]^\l$ such that
$B\not\su \bigcup D$ (we can actually choose $B$ to be disjoint of the
union of $D$)

Using the surjectivity of $\Phi$, fix a graph $G\in \Cal G_\l$ with
$\Phi(G)=B$.

 Suppose to the contrary that $G$ is isomorphic to a subgraph of
$G_\a$ for some $\a<\cont$. Without loss of generality $G$ is a
subgraph of, say, $G_0$. By the covering property of $\Cal F_\a$ we
can find a subset $ F\in [\Cal F_0]^{<\l}$ such that $G\le \bigcup
F$. Because every member of $F$ is an $r$-subgraph of $G_0$, every
member of $F$ is also an $r$-subgraph of (the induced subgraph of
$G_\a$ spanned by) $\bigcup F$. Therefore $\Phi(\bigcup
F)=\bigcup_{A\in F}\phi(A)$ by continuity Claim 2.2 above. Thus 
$$B=\Phi(G)\su \Phi(\bigcup F)=\bigcup_{A\in F}\Phi(A)\su \bigcup D$$

a contradiction to the choice of $B\not\su \bigcup D$.\endproof

\noindent
{\bf Discussion}. Corollary 2.3 has two improvements over
Corollary 1.9: The first is that it extends the result to
singular cardinals. The second is that also it says something
stronger than setting a lower bound for complexity: it says that for
every set of fewer than $\cont$ members in $\Cal G_\mu$ there is a
{\it smaller} graph $G\in \Cal G_\l$ for some regular $\l<\mu$. The
second statement holds for  all regular $\l>\aleph_2$ below the
first fixed point of second order as well, by the  same proof.

\noindent
{\bf Generalization to higher cardinals}

We state now the most general form of representation we know for
$\Cal G$. This involves replacing club guessing using $c_\d$-s of order
type $\om$ with club guessing using $c_\d$-s of order type $\mu$, for
some $\mu$ of cofinality $\aleph_0$. 

Replacing $\om$ and $\Cal P(\om)$ with $\mu$ and $[\mu]^{\aleph_0}$
respectively in the proof of Theorem 1.10, there is
only one difficulty: in establishing condition (d) in the proof, the
argument that $c_\d\su^* E$ implies that
$\phi_{<_1,r}(v,\d)=\phi_{<_2,r}(v,\d)$ breaks down, because a proper
initial segment of $c_\d$ may contain all the members of the subset of
$c_\d$ chosen by $v$ in equation (1) in the proof. 

This difficulty vanishes if we demand that $C_\d\su E$, namely work
with $I=\id^a(\bar C)$ rather than $I=\id^b(\bar C)$. 

Since  normality of $I$ was not used so far (only $\l$-completeness),
all the results so far hold when replacing $I=\id^b(\bar C)$ by the
$\l$-complete $I=\id ^a(\bar C)$. (Normality will be needed in Section
4). 

Thus we have proved:

\ms\noindent{\bf 
2.4  Theorem}:
{\sl
 Suppose that $\l$ is regular. Then for every
$\mu$ satisfying $\mu^+<\l$ there is a surjective homomorphism
$\Phi_\mu:\lng \Cal G_\l,\le_w\rng \to \lng[\mu^{\aleph_0}]^{\le
l},\su\rng$. 
}

\ms\noindent{\bf 
2.5 Corollary}: 
{\sl
For every regular uncountable $\l>\aleph_1$
we have $\cp \Cal G_\l\ge\sup(\{\l^+\} \cup{\mu^{\aleph_0}:\mu^+<\l}\})$. 
}\endproof

Thus, for example, one can have GCH up to $\aleph_\om$ with
$\aleph_\om^{\aleph_0}$ large (see [Ma]), say
$\aleph_{\om+\om+1}$, and in this model the complexity of $\Cal
G_{\aleph_{\om+2}}$ is larger than $\aleph_{\omega+3}$.

\bigbreak
\bigbreak

\noindent
{\bf \S3 Horizontal complexity}

The complexity $\cp \Cal G_\l$ measures the ``depth'' or ``height'' of
the quasi ordering $\lng \Cal G_\l,\le\rng$. Another legitimate way to
measure how complicated a quasi-ordering is, is by estimating its
``width'', namely the supremum of cardinalities of anti-chains. In
this context it means calculating the possible number of pair-wise
non-embeddable graphs in the class.

We use the representation Theorem from Section 2 to prove that in the
class $\Cal G$ under study this number is alway the maximal possible
in every regular $\l>\aleph_1$. This result used no cardinal
arithmetic assumptions. The idea is to use the homomorphism to pull
back antichains from the range. Since the range of $\Phi$ may be too
small (if, say, $\l>\cont$), we have to use $\phi$. The existence of
antichains in the range of $\phi$ follows from the completeness of $I$
for all successor $\l$ and from the normality of $I$ for inaccessible
$\l$, by the recent [GS].  
This recent result of Gitik and Shelah asserts that for all regular
uncountable cardinals $\l>\aleph_1$, every normal ideal concentrating
on the ordinals below $\l$ with cofinality $\aleph_0$ is not $\l^+$
saturated (see [GS]). Not being $\l^+$ saturated means: there are
$\l^+$ positive sets whose pairwise intersections lie in the ideal.
This applies to the guessing ideal $\id^b(\bar C)$ we are using.

\ms\noindent{\bf 
3.1 Theorem}: 
{\sl
For every regular $\l>\aleph_1$ there are
$2^\l$ pair-wise incomparable elements in $\Cal G_\l$ with respect to
$\le_w$. In other words, there are $2^\l$ graphs in $\Cal G_\l$ with
no one of them weakly embeddable in another. 

Furthermore, those graphs can be chosen to be in $\Phi^{-1}(A)$ for a
any given subset $A\in [\Bbb R]^2$.
}

\proof Fix a sequence $\lng S_\a:\a<\l\rng$ of pairwise disjoint
positive sets with respect to $I=\id^b(\ov C)$, $\ov C$ a club guessing
sequence. Let $\Phi$ be the homomorphism from Theorem
1.10 defined using $I$.  

Let $A=\{X,Y\}$ where $X,Y\su
\bar{\Cal P}(\om)$ are two  
distinct, non-empty sets.

Fix a collection $L\su \Cal P(\l)$ of size $2^\l$ with the property
that for every two distinct members $\eta_1,\eta_2\in L$ it holds that
$\eta_1\not\su \eta_2$ and $\eta_2\not\su \eta_1$ ($L$ can be chosen,
for example, to be an independent family over $\l$). For every
$\eta\in L$ we shall find a graph $G_\eta$ with $\Phi(G_\eta)=A$ and
such that none of those graphs is embeddable in another. 

 Suppose $\eta\in L$ is given. Using the surjectivity of $\phi$ find
$G_\eta$ such that $\phi(G)(\d)=X$ iff
$\d\in
\bigcup_{\a\in \eta}S_\a$ and $A_\d=Y$ otherwise. (We neglect
coordinates $\d$ which are not in $S$).

For every $\eta\su\l$ it
holds that $\Phi(G_\eta)=A$, because each $S_\a$ is positive and both
$\eta$ and $\l\sm \eta$ are not empty. Suppose
now that $\eta_1,\eta_2\in L$ are distinct and we
shall show that $G_{\eta_1}\not\le_w G_{\eta_2}$. Since $\phi$ is a
homomorphism, it is enough to show that $\phi(G_{\eta_1})\not\su_I
\phi(G_{\eta_2})$. Suppose then, to the contrary, that for some set
$H\in I$ it holds that $\d\in S\sm H$ implies that
$\phi(G_{\eta_1})(\d)\su\phi(G_{\eta_2})(\d)$. Let $\a\in \eta_1\sm
\eta_2$ be picked. Since $S_\a$ is positive and $H\in I$ the set
$S_\a\sm H$ is positive. For every $\a\in S_\a\sm H$ it holds that
$\phi(G_{\eta_1})(\d)=\{X\}$ and $\phi(G_{\eta_2})=\{Y\}$. Thus
$\phi(G_{\eta_1})\not\su_I\phi(G_{\eta_2})$. \endproof

\ms
\noindent
{\bf Discussion}
The property $\Phi(G)=A$ can be regarded as a strong homogeneity
property: modulo $I$, all rays in neighbourhoods of elements of $G$ converge to
their supremum (when the graphs is well ordered) in exactly one of two
possible convergence rates. Yet, the graphs $G_\eta$ chosen above are
pairwise incomparable.

\bigbreak
\bigbreak
\noindent{\bf \S4 A decomposition Theorem}

Forbidding a few more
countable configurations enables a decomposition theorem: we prove that
such a graph is $r$-indecomposable iff its image is a singleton.

\ms\noindent{\bf 
4.1 Definition}:  Suppose that   $\l>\aleph_1$ is
regular and that $\Phi:\Cal G_\l\to [\Bbb R]^{\le \l}$ is the
homomorphism from Theorem 1.10. Say that   $G\in \Cal
G_\l$ is $(r,\l)$-indecomposable iff for every $\b<\l$ and
$r$-subgraphs $\lng G_\a:\a<\b\rng$ so that $G=\bigcup_{\a<\b}G_\a$,
there is some $\a<\b$ such that $\Phi(G)=\Phi(G_\a)$.

\ms\noindent{\bf 
4.2 Claim}: 
{\sl
Suppose that   $\l>\aleph_1$ is
regular and that $\Phi:\Cal G_\l\to [\Bbb R]^{\le \l}$ is the
homomorphism from Theorem 1.10. If $G\in \Cal G_\l$ and
$|\Phi (G)|=1$ then $G$ is $(r,\l)$-indecomposable.
}

\proof: Suppose $\Phi(G)=\{a\}$.  Fix a positive set $Y\su \l$ and
vertices $v_\d\in G$ for 
$\d\in Y$ so that $\phi(v,\d)=a$. By completeness of $I$ we can find a
fixed $\a<\b$ and a positive $Y'\su Y$ such that $v_\d\in G_\a$ for
all $\d\in Y'$. Since $G_\a$ is an $r$-subgraph, $r(v_\d)\su^* G_\a$
for every $\d\in Y'$. On a measure 1 set equation (8) above holds and tells
us that $\phi(v_\d,d)$ can be computed in either $G_\a$ or in $G$.
Intersecting $Y'$ with this measure 1 set yields a positive measure
$Y''\su Y$ so that $\phi(v_\d,d)=\{a\}$ for $\d\in Y''$. This shows
that $\Phi(G)\su \Phi(G_\a)$. The other inclusion holds because
$G_\a\le G$.\endproof

We may ask if the converse is also true: namely that $\Phi(G)$ can be
represented as a union of fewer than $\l$ $r$-subgraphs $G_\a$,
each with $\phi(G_\a)$ different from, hence properly included in,
$\Phi(g)$ whenever $|\Phi(G)|>1$.  While this is not true for $\Cal
G$, forbidding an additional set of
countable configurations gives a large subclass of $\Cal G$ for which
it is true.

\ms\noindent{\bf 
4.3 Definition}: 
\startitm
\itm Let $\Gamma^*$ be the set of all
countable graphs $G$ in which there is a vertex $v\in G$ with a ray
$R\su G[v]$ and a vertex $u\in R$ with a ray $S\su G[u]$. 
\itm Let $\Cal G^*$ be the subclass of $\Cal G$ resulting by
forbidding all graphs in $\Gamma^*$.

\ms\noindent{\bf 
4.4  Claim}:
{\sl
 If $\l$ and $\Phi$ are as above, then
$\Phi\rest \Cal G^*_\l$ is surjective. Hence the Theorems and
Corollaries proved above for $\Cal G$ hold also for $\Cal G^*$.
}

\proof The graphs constructed in the proof of 1.10 to
demonstrate surjectivity are all in $\Cal G^*$ by condition (d)(i) in
the inductive construction.

\ms\noindent{\bf 
4.5 Claim}:
{\sl
 Suppose $\l$ and $\Phi$ are as above. Suppose
that $G\in \Cal G^*$ and that $|\Phi(G)|>1$. Then there are two
$r$-subgraphs $G_1,G_2$ of $G$ such that $\Phi(G)\not\su\Phi(G_i)$ for
$i=1,2$. 
}

\proof  Let $a,b\in \Phi(G)$ be two distinct members of $\Phi(G)$. For
every $\d<\l$ let $\lng v_{\d,a,i}:i<j(\d,a)\rng$ be an enumeration of
all vertices $v\in G$ with $\phi(v,\d)=a$. Similarly, let $\lng
v_{\d,b,i}:i<j(b,\d)\rng$ enumerate all vertices in $G$ with
$\phi(v,\d)=b$. Let $r$ be a function on $G$ such that $r(v)$ is a ray
in $G[v]$ if such a ray exists and is $\emptyset$ otherwise. Let
$A:=\{v_{\d,a,i}:\d<\l,\,i<j(\d,\a)\}$ and let $G_1=A\cup\{r(v):v\in
A\}$. Since $r(u)=\emptyset$ for every vertex $u\in r(v)$, for $v\in
A$ --- because $\Gamma^*$ is omitted --- we conclude that $G_1$ is an
$r$-subgraphs of $G$. Also, since $r(u)\not=\emptyset$ for every $u\in
B:=\{v_{\d,b,i}:\d<\l,\,i<j(\d,b)\}$, it follows that $G_1\cap
B=\emptyset$. 

Let $C:=G\sm G_1$ and let $G_2=C\cup \{r(v):v\in C\}$. Again,
$r(u)\cap A=\emptyset$ for $u\in C$.

By restricting attention to a measure 1 subset of $\l$ we may assume
that $\phi(v,\d)$ for a vertex $v\in G_1$ (or in $G_2$) remains
unchanged when computed in $G$. Thus we see that $\Phi(G_1)=\{a\}$. On
the other hand, $\phi(v,\d)\not=a$ for all $u\in G_2$, because all
vertices $v\in G$ for which $\phi(v,\d)=a$ for some $\d<\l$ are in
$A$, and $A\cap G_2=\emptyset$. This proves the claim.\endproof

Both claims show that in $\Cal G^*_\l$ a graph is $(r,\l)$-indecomposable if
and only if its image under $\Phi$ is a singleton.

\bigbreak
\bigbreak

\noindent
{\bf References}

\ms\noindent[S] S.~Shelah, {\bf Cardinal Arithmetic}, Oxford University press,
Oxford 1994. 

\ms\noindent[D1] R.~Diestel, {\sl On the problem of finding small subdivision and
homomorphism bases for classes of countable graphs}, Discrete Math. 55
(1985) 21--33.

\ms\noindent[D2]  R.~Diestel, {\sl The end structure of a graph: recent results
and open problems} Discrete  Math. 100 (1992) 313--327

\ms\noindent[DHV] R.~Diestel, R.~Halin and W.~Vogler, {\sl Some remarks on
universal graphs}, Combinatorica 5 (1985) 283--293

\ms\noindent[E] W.~B.~Easton, {\sl Powers of regular cardinals},
Annals of Mathematical logic. 1 (1970) 139--178.

\ms\noindent[FW] M.~Foreman and H.~Woodin, {\sl The generalized
continuum hypothesis can fail everywhere}, Annals of Mathematics 133
(1991) 1--36.

\ms\noindent[GS] M.~Gitik and S.~Shelah, {\sl Less saturated ideals},
preprint, Logic Eprints

\ms\noindent[GrS] R.~Grossberg and S.~Shelah, {\sl On universal
locally finite groups}, Israel Journal of Math. 44 (1983) 289--302

\ms\noindent[HK] A.~Hajnal and P.~Komjath, {\sl Embedding graphs into colored
graphs} Trans. Amer. Math. Soc 307 (1988) 395-409.

\ms\noindent[J] T.~Jech, {\sl Singular cardinal problem: Shelah's
theorem on $2^{\aleph_\om}$}, Bull. London Math. Soc. 24 (1992)
127-139

\ms\noindent[Je] R.~B.~Jensen, {\sl The fine structure of the
constructible universe}, Annals of Math. Logic 4 (1972) 229--308

\ms\noindent[Jo] B.~J\'onsson, {\sl Universal relational systems},
Math. Scand. 4(1956) 193--208

\ms\noindent[k] M~.Kojman, {\sl On universal graphs without cliques or
without large bipartite graphs}, submitted

\ms\noindent[KjS1] M.~kojman and S.~Shelah, {\sl 
nonexistence of universal orders
in many cardinals} Journal of Symb. Logic 57 (1992) 875--891

\ms\noindent[KjS2] M.~kojman and S.~Shelah, {\sl The universality spectrum of
stable unsuprestable  theories} Annals Pure Appl. Logic 58 (1992)
57--72

\ms\noindent[KjS3] M.~kojman and S.~Shelah, {\sl Universal abelian groups} Israel
J. Math, in press.

\ms\noindent[KP] P.~Komjath and J.~Pach, {\sl Universal Graphs without large
bipartite subgraphs}, Mathematika 31 (1984) 282--290

\ms\noindent[KP1]  P.~Komjath and J.~Pach, {\sl Universal elements and the
complexity of certain classes of infinite graphs}, Discrete Math. 95
(1991) 255--270

\ms\noindent[KS] P.~Komjath and S.~Shelah, {\sl Universal graphs without large
cliques}, preprint

\ms\noindent[L] R.~Laver, {\sl On Fraisse's order type conjecture},
Annals of Mathematics 93 (1971) 89--111

\ms\noindent[M] A.~Macintyre, {\sl Existentially closed structures and
Jensen's principle $\diamondsuit$}, Israel Journal of Math. 25 (1976) 202--210

\ms\noindent[Ma] M.~Magidor, {\sl On the singular cardinal problem II},
Annals of Mathematics 106 (1977), 517--547

\ms\noindent[Me] Alan.~H.~Mekler, {\sl Universal structures in power
$\aleph_1$}, Journal of Symbolic Logic 55(2) (1990) 466--477.

\ms\noindent[MV] M.~Morley and  R.~L.~Vaught, {\sl Homogeneous
universal models}, Math. Scand. 11 (1962) 35--57

\ms\noindent[R] R.~Rado, {\sl Universal graphs and universal functions}, Acta.
Arith. 9. (1964) 331-340

\ms\noindent[S1] S.~Shelah, {\sl Independence results}, Journal of
Symbolic Logic 45 (3) (1980) 563--573

\ms\noindent[S2] S.~Shelah, {\sl On universal graphs without instances of CH},
Annals of Pure. Appl. Logic 26 (1984) 75--87.

\ms\noindent[S3] S.~Shelah, {\sl Universal graphs without instances of CH:
revisitied}, Israel J. Math 70 (1990) 69--81

\ms\noindent[S4] S.~Shelah, {\sl Universal in $<\l$-stable abelian
groups}, preprint

\ms \noindent[S5] S.~Shelah, {\sl Non-existence of universals for
classes like reduced torsion free abelian groups under non necessarily
pure embeddings}, in preparation
\end